\documentclass[11pt]{amsart}

\usepackage{amsmath}
\usepackage{amssymb}
\usepackage{graphicx}

\theoremstyle{definition}

\numberwithin{equation}{section}

\title[Solvability for operators with Laplacian and
bi-Laplacian]{Solvability in the sense of sequences for certain non-
Fredholm operators with a drift and Laplace and bi-Laplace operators}

\author[V. Vougalter]{Vitali Vougalter}

\address[V. Vougalter]{Department of Mathematics,  University of Toronto,
Toronto, Ontario, M5S 2E4, Canada}
\email{{\tt vitali@math.toronto.edu}}



\begin{document}

\begin{abstract}
 We study the solvability of some linear
nonhomogeneous elliptic problems and establish that under certain technical
assumptions the convergence in $L^{2}$ of their right-hand sides yields the
existence and the convergence in $H^{4}$ of the solutions. The equations
contain fourth order differential operators with or without the Fredholm
property, in particular the second and the fourth derivative operators,
on the whole real line or on a finite interval with periodic
boundary conditions. We establish that the transport term involved in these
problems provides the regularization of the solutions.
\end{abstract}

\maketitle


\setcounter{equation}{0}

\section{Introduction}

Let us consider the problem
\begin{equation}
\label{eq1}
(-\Delta + V(x))u - a u=f
\end{equation}
with $u \in E= H^{2}({\mathbb R}^{d})$ and  $f \in F=
L^{2}({\mathbb R}^{d}), \ d\in {\mathbb N}$. Here $a$ is a constant, and the function
$V(x)$ tends to $0$ at the infinity. 
When $a \geq 0$, the origin belongs to the essential spectrum of the operator $A : E \to F$,
which corresponds to the left side of equation (\ref{eq1}).
As a consequence, this operator fails to 
satisfy the Fredholm property. Its image is not closed, for $d>1$
the dimension of its kernel and the codimension of its image are
not finite. The present article deals with the studies of the certain
properties of the operators of this kind.
Note that the elliptic equations involving the non-Fredholm operators were
treated extensively in recent years (see ~\cite{EV20}, ~\cite{EV21},
~\cite{EV211}, ~\cite{EV22}, ~\cite{V2011}, ~\cite{VKMP02}, ~\cite{VV131}, ~\cite{VV220}, ~\cite{V25}, ~\cite{V26}, ~\cite{VV10}, 
~\cite{VV19}, ~\cite{VV21}, ~\cite{VV22}, also ~\cite{B88})
along with their
potential applications to the theory of reaction-diffusion
problems (see ~\cite{DMV05}, ~\cite{DMV08}). Fredholm structures, topological
invariants and their applications were discussed in ~\cite{E09}.
The articles ~\cite{GS05} and ~\cite{RS01} are important for the understanding of the
Fredholm
and properness properties of the quasilinear elliptic systems of the second
order and of the operators of this kind on ${\mathbb R}^{N}$.
The work  ~\cite{GS10} deals with the exponential decay and Fredholm properties in the second-order
quasilinear elliptic systems of equations.
In particular, when $a$ is trivial, the operator $A$ satisfies the Fredholm property
in the certain properly chosen weighted spaces (see \cite{Amrouche1997},
\cite{Amrouche2008}, \cite{B88}, \cite{Bolley1993}, \cite{Bolley2001}).
However, the situation when $a$ is nontrivial is considerably
different and the method developed in these articles cannot be
used.

One of the major questions concerning the problems involving non-Fredholm operators
is their solvability. Let us address it in the following
 setting. Suppose $f_{n}$ is a sequence of
 functions in the image of the operator $A$, such that
$f_{n}\to f$ in $L^{2}({\mathbb R}^{d})$ as $n\to \infty$. We denote
 by $u_n$  a sequence of functions from $H^{2}({\mathbb R}^{d})$, so that
$$
Au_{n}=f_{n}, \ n\in {\mathbb N}.
$$
Since the operator $A$ does not satisfy the Fredholm property,
 the sequence $u_n$ may not be convergent.
 Let us call a sequence $u_n$, such that $A u_n \to f$ in
 $L^{2}({\mathbb R}^{d})$ a solution in
 the sense of sequences of the equation $Au=f$ (see \cite{V2011}). If this
 sequence tends to a function $u_0$ in the norm of the space
 $E$, then $u_0$ is a solution of this problem. The solution in
 the sense of sequences is equivalent in this sense to the
 usual solution. However, in the situation of non-Fredholm
operators this convergence may not hold or it can occur in a certain
 weaker sense. In this case, the solution in the sense of sequences may
 not imply the existence of the usual solution. In the present work we
 will find the sufficient conditions of equivalence of solutions in the
 sense of sequences and the usual solutions. In the other
 words, we will determine the conditions on sequences $f_n$ under which the
 corresponding sequences $u_n$ are strongly convergent.

In the first part of the work we consider the problem with the drift
term
\begin{equation}
\label{eq2}
-\frac{d^{2}u}{dx^{2}}+\frac{d^{4}u}{dx^{4}}-b\frac{du}{dx}-au=f(x), \quad x\in {\mathbb R},
\end{equation}
where $a\geq 0$ and $b\in {\mathbb R}, \ b\neq 0$ are the constants and the
right side is square integrable. 
The equation with the transport in the context of the
Darcy's law describing the fluid motion in the porous medium was studied
in ~\cite{VV10}. The drift term is significant for the understanding of the emergence and
propagation of patterns arising in the theory of speciation (see ~\cite{VV13}).
Nonlinear propagation phenomena for the reaction-diffusion type problems
involving the transport term were discussed in ~\cite{BHN05}.
Existence of solutions for certain non-Fredholm integro-differential
equations with the bi-Laplacian was established in ~\cite{VV21}. The work
~\cite{EV21} is devoted to the solvability in the sense of sequences for some
fourth order non-Fredholm operators. 
Evidently, the operator contained in the left side of
(\ref{eq2})
\begin{equation}
\label{lab}
L_{a, \ b}:=-\frac{d^{2}}{dx^{2}}+\frac{d^{4}}{dx^{4}}-b\frac{d}{dx}-a: \quad H^{4}({\mathbb R})
\to L^{2}({\mathbb R})
\end{equation}
is non-selfadjoint.
By virtue of the standard Fourier transform
\begin{equation}
\label{ft}
\widehat{f}(p):=\frac{1}{\sqrt{2\pi}}\int_{-\infty}^{\infty}f(x)
e^{-ipx}dx, \quad p\in {\mathbb R}
\end{equation}
it can be trivially checked that the essential spectrum of the operator $L_{a, \ b}$
is given by
$$
\lambda_{a, \ b}(p):=p^{2}+p^{4}-a-ibp, \quad p\in {\mathbb R}.
$$
Clearly, for $a>0$ the operator $L_{a, \ b}$ is Fredholm since its essential spectrum does not contain the origin. But if $a=0$, our operator
$L_{0, \ b}$ fails to satisfy the Fredholm property because the origin belongs to  its essential spectrum.

Obviously, in the absence of the drift term we are dealing with the
self-adjoint operator
$$
-\frac{d^{2}}{dx^{2}}+\frac{d^{4}}{dx^{4}}-a: \quad H^{4}({\mathbb R})\to L^{2}({\mathbb R}),
\quad a\geq 0,
$$
which is non-Fredholm. Thus, the
introduction of the transport term provides the regularization for the solutions
of our equation.
We write down the corresponding sequence of approximate equations with
$m\in {\mathbb N}$, namely
\begin{equation}
\label{eq3}
-\frac{d^{2}u_{m}}{dx^{2}}+\frac{d^{4}u_{m}}{dx^{4}}-b\frac{du_{m}}{dx}-au_{m}=f_{m}(x), \quad x\in
{\mathbb R},
\end{equation}
where $a\geq 0$ and $b\in {\mathbb R}, \ b\neq 0$ are the constants.
The right sides of (\ref{eq3}) converge to the right side of (\ref{eq2}) in
$L^{2}({\mathbb R})$ as $m\to \infty$. Let us introduce the inner product of two
functions
\begin{equation}
\label{ip}
(f(x),g(x))_{L^{2}({\mathbb R})}:=\int_{-\infty}^{\infty}f(x)\bar{g}(x)dx,
\end{equation}
with a slight abuse of notations when these functions are not square integrable
on the real line. Indeed, if $f(x)\in L^{1}({\mathbb R})$ and
$g(x)$ is bounded, then the integral in the right side of (\ref{ip}) makes sense,
like for instance in the situation of the functions involved in the
orthogonality conditions (\ref{or1}) and (\ref{or2}) of Theorems 1.1 and 1.2
below.
For our problem (\ref{eq2}) on the finite interval $I:=[0, 2\pi]$ with periodic
boundary conditions (see (\ref{h1})), we will use the inner product analogous
to (\ref{ip}), replacing the real line with $I$.
In the first part of the present work we will consider the space
$H^{4}({\mathbb R})$ equipped with the norm
\begin{equation}
\label{h2n}
\|u\|_{H^{4}({\mathbb R})}^{2}:=\|u\|_{L^{2}({\mathbb R})}^{2}+ \Bigg\|
\frac{d^{4}u}{dx^{4}}\Bigg\|_{L^{2}({\mathbb R})}^{2}.
\end{equation}
When working with the norm $H^{4}(I)$ further down, we will replace ${\mathbb R}$
with $I$ in formula (\ref{h2n}).
Our first main proposition is as follows.

\bigskip

\noindent
{\bf Theorem 1.1.} {\it Suppose the constants
$a\geq 0, \ b\in {\mathbb R}, \ b\neq 0$ and $f(x)\in L^{2}({\mathbb R})$.

\medskip

\noindent
a) When $a>0$, problem (\ref{eq2}) has a unique solution
$u(x)\in H^{4}({\mathbb R})$.

\medskip

\noindent
b) When $a=0$, let in addition $xf(x)\in L^{1}({\mathbb R})$. Then equation
(\ref{eq2}) admits a unique solution
$u(x)\in H^{4}({\mathbb R})$ if and only if the orthogonality relation
\begin{equation}
\label{or1}
(f(x), 1)_{L^{2}({\mathbb R})}=0
\end{equation}
holds.}

\bigskip

Note that the quantity in the left side of (\ref{or1}) is well defined by means of
the result of Lemma 4.1 of ~\cite{VV19}.
Evidently, the argument of the case a) of the theorem above is not based on
the orthogonality relations.
Our next result is about the solvability in the sense of
sequences for our problem on the whole real line.

\bigskip

\noindent
{\bf Theorem 1.2.} {\it Let the constants
$a\geq 0, \ b\in {\mathbb R}, \ b\neq 0$ and
$m\in {\mathbb N}, \ f_{m}(x)\in L^{2}({\mathbb R})$, such that
$f_{m}(x)\to f(x)$ in $L^{2}({\mathbb R})$ as $m\to \infty$.

\medskip

\noindent
a) For $a>0$, problems (\ref{eq2}) and (\ref{eq3}) possess unique solutions
$u(x)\in H^{4}({\mathbb R})$ and $u_{m}(x)\in H^{4}({\mathbb R})$ respectively,
so that $u_{m}(x)\to u(x)$ in $H^{4}({\mathbb R})$ as $m\to \infty$.

\medskip

\noindent
b) For $a=0$, let $xf_{m}(x)\in L^{1}({\mathbb R})$, such that 
$xf_{m}(x)\to xf(x)$ in $L^{1}({\mathbb R})$ as $m\to \infty$. Moreover,
\begin{equation}
\label{or2}
(f_{m}(x), 1)_{L^{2}({\mathbb R})}=0, \quad m\in {\mathbb N}
\end{equation}
holds. Then equations (\ref{eq2}) and (\ref{eq3}) have unique solutions
$u(x)\in H^{4}({\mathbb R})$ and $u_{m}(x)\in H^{4}({\mathbb R})$ respectively,
so that $u_{m}(x)\to u(x)$ in $H^{4}({\mathbb R})$ as $m\to \infty$.}

\bigskip

The second part of the article deals with the studies of our problem on the
finite interval with the periodic boundary conditions (see (\ref{h1})), i.e.
$I=[0, \ 2\pi]$, namely
\begin{equation}
\label{eq4}
-\frac{d^{2}u}{dx^{2}}+\frac{d^{4}u}{dx^{4}}-b\frac{du}{dx}-au=f(x), \quad x\in I,
\end{equation}
where $a\geq 0$ and $b\in {\mathbb R}, \ b\neq 0$ are the constants and the
right side of (\ref{eq4}) is continuous and periodic. Evidently,
\begin{equation}
\label{f12}
\|f\|_{L^{1}(I)}\leq 2\pi \|f\|_{C(I)}<\infty, \quad
\|f\|_{L^{2}(I)}\leq \sqrt{2\pi} \|f\|_{C(I)}<\infty.
\end{equation}
Hence, $f(x)\in L^{2}(I)$ as well. Let us use the Fourier
transform
\begin{equation}
\label{fti}
f_{n}:=\frac{1}{\sqrt{2\pi}}\int_{0}^{2\pi}f(x)e^{-inx}dx, \quad n\in {\mathbb Z},
\end{equation}
such that
$$
f(x)=\sum_{n=-\infty}^{\infty}f_{n}\frac{e^{inx}}{\sqrt{2\pi}}.
$$
Clearly, the non-selfadjoint operator involved in the left side of
(\ref{eq4})
\begin{equation}
\label{labi}
l_{a, \ b}:=-\frac{d^{2}}{dx^{2}}+\frac{d^{4}}{dx^{4}}-b\frac{d}{dx}-a: \quad H^{4}(I)\to
L^{2}(I)
\end{equation}
is Fredholm. By virtue of (\ref{fti}), it can be easily verified that
the spectrum of $l_{a, \ b}$ is given by
$$
\lambda_{a, \ b}(n):=n^{2}+n^{4}-a-ibn, \quad n\in {\mathbb Z}.
$$
The corresponding eigenfunctions are the standard Fourier harmonics
$$
\frac{e^{inx}}{\sqrt{2\pi}}, \ n\in {\mathbb Z}.
$$
Note that the eigenvalues of the operator  $l_{a, \ b}$ are simple.
The corresponding function space here $H^{4}(I)$ is 
$$  
\{u(x):I\to {\mathbb C} \ | \ u(x), u''''(x)\in L^{2}(I), \quad
u(0)=u(2\pi), \quad u'(0)=u'(2\pi),
$$
\begin{equation}
\label{h1} 
u''(0)=u''(2\pi), \quad u'''(0)=u'''(2\pi)\}.
\end{equation}
For the technical goals, let us introduce the following auxiliary constrained
subspace
\begin{equation}
\label{H0}
H_{0}^{4}(I)=\{u(x)\in H^{4}(I) \ | \ (u(x),1)_{L^{2}(I)}=0 \},
\end{equation}
which is a Hilbert spaces as well (see e.g. Chapter 2.1 of ~\cite{HS96}).
Obviously, if $a>0$, the kernel of the operator $l_{a, \ b}$  is trivial.
When $a=0$, we consider
$$
l_{0, \ b}: \quad H_{0}^{4}(I)\to L^{2}(I).
$$
This operator also has a trivial kernel.
Let us write down the corresponding sequence of the approximate equations with
$m\in {\mathbb N}$, namely
\begin{equation}
\label{eq5}
-\frac{d^{2}u_{m}}{dx^{2}}+\frac{d^{4}u_{m}}{dx^{4}}-b\frac{du_{m}}{dx}-au_{m}=f_{m}(x), \quad x\in I,
\end{equation}
where $a\geq 0, \ b\in {\mathbb R}, \ b\neq 0$ are the constants.
The right sides of (\ref{eq5}) are continuous, periodic and tend to the
right side of (\ref{eq4}) in $C(I)$ as $m\to \infty$. The purpose of Theorems 1.3
and 1.4 below is to show the formal similarity of the results on the
finite interval with periodic boundary conditions to the ones derived for the situation
on the whole real in Theorems 1.1 and 1.2 above.

\bigskip

\noindent
{\bf Theorem 1.3.} {\it Let the constants
$a\geq 0, \ b\in {\mathbb R}, \ b\neq 0$ and
$f(0)=f(2\pi), \ f(x)\in C(I)$.

\medskip

\noindent
a) When $a>0$, problem (\ref{eq4}) has a unique solution
$u(x)\in H^{4}(I)$.

\medskip

\noindent
b) When $a=0$, equation (\ref{eq4}) admits a unique solution
$u(x)\in H_{0}^{4}(I)$ if and only if the orthogonality relation
\begin{equation}
\label{or3}
(f(x), 1)_{L^{2}(I)}=0
\end{equation}
holds.}

\bigskip

The final main statement of the work deals with the solvability in the
sense of sequences for our problem on the finite interval $I$.

\bigskip

\noindent
{\bf Theorem 1.4.} {\it Let the constants
$a\geq 0, \ b\in {\mathbb R}, \ b\neq 0$ and $m\in{\mathbb N}$,
such that $f_{m}(0)=f_{m}(2\pi)$. Moreover, $f_{m}(x)\in C(I)$ and
$f_{m}(x)\to f(x)$ in $C(I)$ as $m\to \infty$.

\medskip

\noindent
a) For $a>0$ equations (\ref{eq4}) and (\ref{eq5}) possess unique solutions
$u(x)\in H^{4}(I)$ and $u_{m}(x)\in H^{4}(I)$ respectively, such that
$u_{m}(x)\to u(x)$ in $H^{4}(I)$ as $m\to \infty$.

\medskip

\noindent
b) For $a=0$, let
\begin{equation}
\label{or4}
(f_{m}(x), 1)_{L^{2}(I)}=0, \quad m\in {\mathbb N}.
\end{equation}
Then problems (\ref{eq4}) and (\ref{eq5}) have unique solutions
$u(x)\in H_{0}^{4}(I)$ and $u_{m}(x)\in H_{0}^{4}(I)$ respectively, such that
$u_{m}(x)\to u(x)$ in $H_{0}^{4}(I)$ as $m\to \infty$.}

\bigskip

Let us note that in the situations a) of Theorems 1.3 and 1.4 above the argument is not based
on the orthogonality conditions.
Clearly, when there are no drift terms in our equations, the
case is more singular (see formulas
(\ref{unfn}) and (\ref{ufmn}) below with
$a=n_{0}^{2}+n_{0}^{4}, \ n_{0}\in {\mathbb N}$).

\bigskip


\setcounter{equation}{0}

\section{The whole real line case}

\noindent
{\it Proof of Theorem 1.1.} First we establish that it would be
sufficient to solve our problem in $L^{2}({\mathbb R})$. Suppose $u(x)$ is
a square integrable solution of (\ref{eq2}) on the whole real line. Directly
from this equation under the stated assumptions we obtain that
$$
-\frac{d^{2}u}{dx^{2}}+\frac{d^{4}u}{dx^{4}}-b\frac{du}{dx}\in L^{2}({\mathbb R})
$$
as well. Let us apply here the standard Fourier transform (\ref{ft}). Hence,
$(p^{2}+p^{4}-ibp){\widehat u}(p)\in L^{2}({\mathbb R})$, such that
$\displaystyle{\int_{-\infty}^{\infty}p^{8}|{\widehat u}(p)|^{2}dp<\infty}$.
Thus, $\displaystyle{\frac{d^{4}u}{dx^{4}}\in L^{2}({\mathbb R})}$.
This means that $u(x)\in H^{4}({\mathbb R})$ as well.

Let us demonstrate the uniqueness of solutions for our problem. Suppose
that
$$
u_{1}(x), \ u_{2}(x)\in H^{4}({\mathbb R})
$$
solve equation (\ref{eq2}).
Clearly, their difference $w(x):=u_{1}(x)-u_{2}(x)\in H^{4}({\mathbb R})$.
It is a solution of the homogeneous equation
$$
-\frac{d^{2}w}{dx^{2}}+\frac{d^{4}w}{dx^{4}}-b\frac{dw}{dx}-aw=0.
$$
Evidently, the operator $L_{a, \ b}: H^{4}({\mathbb R})\to L^{2}({\mathbb R})$ defined in (\ref{lab}) does not have
any nontrivial zero modes. Therefore, the function $w(x)\equiv 0$ on ${\mathbb R}$.

Apply the standard Fourier transform (\ref{ft}) to both sides of problem
(\ref{eq2}). This gives us
\begin{equation}
\label{upfp}
{\widehat u}(p)=\frac{{\widehat f}(p)}{p^{2}+p^{4}-a-ibp}.
\end{equation}
Hence,
\begin{equation}
\label{upfpn}
\|u\|_{L^{2}({\mathbb R})}^{2}=\int_{-\infty}^{\infty}\frac{|{\widehat f}(p)|^{2}}
{(p^{2}+p^{4}-a)^{2}+b^{2}p^{2}}dp.
\end{equation}
Consider first the situation a) of the theorem. By means of (\ref{upfpn}), we obtain
$$
\|u\|_{L^{2}({\mathbb R})}^{2}\leq \frac{1}{C}\|f\|_{L^{2}({\mathbb R})}^{2}<\infty
$$
as we assume. Here and below $C$ will stand for
a finite, positive constant.

Let us turn our attention to the case when the parameter $a$ is trivial.
By means of (\ref{upfp}), we have
\begin{equation}
\label{upfp0}
{\widehat u}(p)=\frac{{\widehat f}(p)}{p^{2}+p^{4}-ibp}\chi_{\{|p|\leq 1\}}+
\frac{{\widehat f}(p)}{p^{2}+p^{4}-ibp}\chi_{\{|p|>1\}}.
\end{equation}
Here and further down $\chi_{A}$ will designate the characteristic function
of a set $A\subseteq {\mathbb R}$.
Obviously, the second term in the right side of (\ref{upfp0}) can be bounded
from above in the absolute value by
$\displaystyle{\frac{|{\widehat f}(p)|}{|b|}\in L^{2}({\mathbb R})}$ since
$f(x)$ is square integrable on ${\mathbb R}$ as assumed.
Note that
$$
{\widehat f}(p)={\widehat f}(0)+\int_{0}^{p}\frac{d{\widehat f}(s)}{ds}ds.
$$
Thus, the first term in the right side of (\ref{upfp0}) is given by
\begin{equation}
\label{fp0}
\frac{{\widehat f}(0)}{p^{2}+p^{4}-ibp}\chi_{\{|p|\leq 1\}}+\frac{\int_{0}^{p}
\frac{d{\widehat f}(s)}{ds}ds}{p^{2}+p^{4}-ibp}\chi_{\{|p|\leq 1\}}.
\end{equation}
Using definition (\ref{ft}) of the standard Fourier transform, we easily
derive that
$$
\Bigg|\frac{d{\widehat f}(p)}{dp}\Bigg|\leq \frac{1}{\sqrt{2\pi}}\|xf(x)\|_
{L^{1}({\mathbb R})}.
$$
Hence, the second term in (\ref{fp0}) can be estimated from above in the
absolute value by
$$
\frac{1}{\sqrt{2\pi}}\frac{\|xf(x)\|_{L^{1}({\mathbb R})}}{|b|}
\chi_{\{|p|\leq 1\}}\in L^{2}({\mathbb R}).
$$
It can be easily verified that the first term in (\ref{fp0}) belongs to $L^{2}({\mathbb R})$ if
and only if ${\widehat f}(0)$ vanishes. This is equivalent to (\ref{or1}).  \hspace{4cm} $\Box$

\bigskip

We turn our attention to establishing the solvability in the sense of
sequences for our equation on ${\mathbb R}$.

\bigskip

\noindent
{\it Proof of Theorem 1.2.} Let us first suppose that problems (\ref{eq2})
and (\ref{eq3}) have unique solutions $u(x)\in H^{4}({\mathbb R})$ and
$u_{m}(x)\in H^{4}({\mathbb R}), \ m\in {\mathbb N}$ respectively, such that
$u_{m}(x)\to u(x)$ in $L^{2}({\mathbb R})$ as $m\to \infty$. This will yield
that $u_{m}(x)$ also tends to $u(x)$ in $H^{4}({\mathbb R})$ as
$m\to \infty$. Using (\ref{eq2}) and (\ref{eq3}), we easily derive that
$$
\Bigg\|-\frac{d^{2}}{dx^{2}}(u_{m}-u)+\frac{d^{4}}{dx^{4}}(u_{m}-u)-b\frac{d(u_{m}-u)}{dx}\Bigg\|_
{L^{2}({\mathbb R})}\leq
$$
\begin{equation}
\label{d4dx4umu}  
\|f_{m}-f\|_{L^{2}({\mathbb R})}+a\|u_{m}-u\|_{L^{2}({\mathbb R})}.
\end{equation}
The right side of (\ref{d4dx4umu}) converges to zero as $m\to \infty$ as
we assume. By virtue of the standard Fourier transform (\ref{ft}), we
easily obtain that
$$
\int_{-\infty}^{\infty}p^{8}|{\widehat u_{m}}(p)-{\widehat u}(p)|^{2}dp\to 0, \quad
m\to \infty.
$$
Hence,
$\displaystyle{\frac{d^{4}u_{m}}{dx^{4}}\to\frac{d^{4}u}{dx^{4}}}$ in
$L^{2}({\mathbb R})$ as $m\to \infty$. This means that
$u_{m}(x)\to u(x)$ in $H^{4}({\mathbb R})$ as $m\to \infty$ as well.

Let us apply the standard Fourier transform (\ref{ft}) to both sides of
(\ref{eq3}). This yields
\begin{equation}
\label{upfpm}
{\widehat u_{m}}(p)=\frac{{\widehat f_{m}}(p)}{p^{2}+p^{4}-a-ibp}, \quad m\in
{\mathbb N}.
\end{equation}
First we consider the situation  a) of the theorem. Recall the part a) of
Theorem 1.1.  Thus, problems (\ref{eq2}) and (\ref{eq3}) possess unique solutions
$u(x)\in H^{4}({\mathbb R})$ and
$u_{m}(x)\in H^{4}({\mathbb R}), \ m\in {\mathbb N}$ respectively.
By means of (\ref{upfpm})
along with (\ref{upfp}), we derive that
$$
\|u_{m}-u\|_{L^{2}({\mathbb R})}^{2}=\int_{-\infty}^{\infty}
\frac{|{\widehat f_{m}}(p)-{\widehat f}(p)|^{2}}{(p^{2}+p^{4}-a)^{2}+b^{2}p^{2}}dp,
$$
such that
$$
\|u_{m}-u\|_{L^{2}({\mathbb R})}\leq \frac{1}{C}\|f_{m}-f\|_{L^{2}({\mathbb R})}\to 0, \quad
m\to \infty
$$
as assumed.
Thus, in the case when $a>0$ we have $u_{m}(x)\to u(x)$ in
$H^{4}({\mathbb R})$ as $m\to \infty$ by virtue of the reasoning above.

We conclude the proof of our theorem by discussing the situation when
the constant $a$ vanishes. Recall the result of the part a) of
Lemma 3.3 of ~\cite{VV131}. Hence, under the given conditions
\begin{equation}
\label{orl}
(f(x), 1)_{L^{2}({\mathbb R})}=0
\end{equation}
holds. Then by means of the part b) of Theorem 1.1, equations (\ref{eq2}) and
(\ref{eq3}) admit unique solutions $u(x)\in H^{4}({\mathbb R})$ and
$u_{m}(x)\in H^{4}({\mathbb R}), \ m\in {\mathbb N}$ respectively if $a$ is trivial.
It easily follows from (\ref{upfpm}) and (\ref{upfp}) that
\begin{equation}
\label{umun}
{\widehat u_{m}}(p)-{\widehat u}(p)=\frac{{\widehat f_{m}}(p)-{\widehat f}(p)}
{p^{2}+p^{4}-ibp}\chi_{\{|p|\leq 1\}}+\frac{{\widehat f_{m}}(p)-{\widehat f}(p)}
{p^{2}+p^{4}-ibp}\chi_{\{|p|>1\}}.
\end{equation}
Evidently, the second term in the right side of (\ref{umun}) can be estimated
from above in the absolute value by
$\displaystyle{\frac{|{\widehat f_{m}}(p)-{\widehat f}(p)|}{|b|}}$. Then
$$
\Bigg\|\frac{{\widehat f_{m}}(p)-{\widehat f}(p)}{p^{2}+p^{4}-ibp}\chi_{\{|p|>1\}}
\Bigg\|_{L^{2}({\mathbb R})}\leq 
\frac{1}{|b|}\|f_{m}-f\|_{L^{2}({\mathbb R})}\to 0, \quad m\to \infty
$$
as we assume. Orthogonality relations (\ref{orl}) and (\ref{or2}) imply that
$$
{\widehat f}(0)=0, \quad {\widehat f_{m}}(0)=0, \quad m\in {\mathbb N}.
$$
Clearly,
\begin{equation}
\label{fmf}
{\widehat f}(p)=\int_{0}^{p}\frac{d{\widehat f}(s)}{ds}ds, \quad
{\widehat f_{m}}(p)=\int_{0}^{p}\frac{d{\widehat f_{m}}(s)}{ds}ds, \quad
m\in {\mathbb N}.
\end{equation}
Therefore, it remains to analyze the term
$$
\frac{\int_{0}^{p}[\frac{d{\widehat f_{m}}(s)}{ds}-\frac{d{\widehat f}(s)}{ds}]
ds}{p^{2}+p^{4}-ibp}\chi_{\{|p|\leq 1\}}.
$$
By virtue of the definition of the standard Fourier transform (\ref{ft}), we
easily derive that
$$
\Big|\frac{d{\widehat f_{m}}(p)}{dp}-\frac{d{\widehat f}(p)}{dp}\Big|\leq
\frac{1}{\sqrt{2\pi}}\|xf_{m}(x)-xf(x)\|_{L^{1}({\mathbb R})}.
$$
Thus,
$$
\Bigg|\frac{\int_{0}^{p}[\frac{d{\widehat f_{m}}(s)}{ds}-\frac{d{\widehat f}(s)}
{ds}]ds}{p^{2}+p^{4}-ibp}\chi_{\{|p|\leq 1\}}\Bigg|\leq
\frac{\|xf_{m}(x)-xf(x)\|_{L^{1}({\mathbb R})}}{\sqrt{2\pi}|b|}\chi_{\{|p|\leq 1\}}.
$$
Obviously,
$$
\Bigg\|\frac{\int_{0}^{p}[\frac{d{\widehat f_{m}}(s)}{ds}-\frac{d{\widehat f}(s)}
{ds}]ds}{p^{2}+p^{4}-ibp}\chi_{\{|p|\leq 1\}}\Bigg\|_{L^{2}({\mathbb R})}\leq
\frac{\|xf_{m}(x)-xf(x)\|_{L^{1}({\mathbb R})}}{\sqrt{\pi}|b|}\to 0
$$
as $m\to \infty$ due to the one of the given conditions. This means that
$u_{m}(x)\to u(x)$ in
$L^{2}({\mathbb R})$ as $m\to\infty$. According to the argument above, we
have $u_{m}(x)\to u(x)$ in $H^{4}({\mathbb R})$ as
$m\to\infty$ in the case b) of our theorem.  \hspace{1cm}   $\Box$

\bigskip


\setcounter{equation}{0}

\section{The equation on the finite interval}

\noindent
{\it Proof of Theorem 1.3.} First we demonstrate that it would be
sufficient to solve our equation in $L^{2}(I)$. Clearly, if $u(x)$ is a square
integrable solution of (\ref{eq4}), periodic on $I$ along with its derivatives
up to the third order inclusively,
directly from our problem under the stated assumptions we obtain that
$$
-\frac{d^{2}u}{dx^{2}}+\frac{d^{4}u}{dx^{4}}-b\frac{du}{dx}\in L^{2}(I).
$$
Let us use (\ref{fti}).  Thus,
$(n^{2}+n^{4}-ibn)u_{n}\in l^{2}$, which implies that
$$
\displaystyle{\sum_{n=-\infty}^{\infty}n^{8}|u_{n}|^{2}<\infty}.
$$
Therefore, $\displaystyle{\frac{d^{4}u}{dx^{4}}\in L^{2}(I)}$,
such that $u(x)\in H^{4}(I)$ as well.

Let us establish the uniqueness of solutions of (\ref{eq4}). We consider the
case when $a>0$. If $a$ is trivial, we are able to apply the similar ideas in the
constrained subspace $H_{0}^{4}(I)$. Suppose
that $u_{1}(x), \ u_{2}(x)\in H^{4}(I)$ satisfy (\ref{eq4}).
Then their difference $w(x):=u_{1}(x)-u_{2}(x)\in H^{4}(I)$. Evidently, it solves
the homogeneous problem
$$
-\frac{d^{2}w}{dx^{2}}+\frac{d^{4}w}{dx^{4}}-b\frac{dw}{dx}-aw=0.
$$
Note that the operator $l_{a, \ b}: H^{4}(I)\to L^{2}(I)$ introduced in
(\ref{labi}) does not possess any nontrivial zero modes. This means that the
function $w(x)$ vanishes identically in $I$.

Let us apply the Fourier transform (\ref{fti}) to both sides of equation
(\ref{eq4}). We arrive at
\begin{equation}
\label{unfn}
u_{n}=\frac{f_{n}}{n^{2}+n^{4}-a-ibn}, \quad n\in {\mathbb Z}.
\end{equation}
Hence,
\begin{equation}
\label{unfnn}
\|u\|_{L^{2}(I)}^{2}=\sum_{n=-\infty}^{\infty}\frac{|f_{n}|^{2}}{(n^{2}+n^{4}-a)^{2}+b^{2}n^{2}}.
\end{equation}
First we treat the case a) of our theorem. It easily follows from (\ref{unfnn})
that
$$
\|u\|_{L^{2}(I)}^{2}\leq \frac{1}{C}\|f\|_{L^{2}(I)}^{2}<\infty
$$
as assumed (see (\ref{f12})). By virtue of the argument
above, we have $u(x)\in H^{4}(I)$.

We conclude the proof of our theorem by covering the case when $a=0$.
From (\ref{unfn}) we easily deduce that
\begin{equation}
\label{unfn0}
u_{n}=\frac{f_{n}}{n^{2}+n^{4}-ibn}, \quad n\in {\mathbb Z}.
\end{equation}
Then
\begin{equation}
\label{unfnn0}
\|u\|_{L^{2}(I)}^{2}=\sum_{n=-\infty}^{\infty}\frac{|f_{n}|^{2}}{n^{2}[n^{2}(1+n^{2})^{2}+b^{2}]}.
\end{equation}
Suppose that $u(x)\in H_{0}^{4}(I)$.  This yields
\begin{equation}
\label{f0}
f_{0}=0.
\end{equation}
On the other hand, if (\ref{f0}) holds, we easily obtain
$$
\|u\|_{L^{2}(I)}^{2}=\sum_{n\in {\mathbb Z}, \ n\neq 0}\frac{|f_{n}|^{2}}
{n^{2}[n^{2}(1+n^{2})^{2}+b^{2}]}\leq \frac{1}{b^{2}}\|f\|_{L^{2}(I)}^{2}<\infty
$$
via the one of the given conditions along with (\ref{f12}). The reasoning above
gives us that $u(x)\in H_{0}^{4}(I)$. Obviously, (\ref{f0}) is equivalent
to orthogonality relation (\ref{or3}).  \hspace{9.5cm} $\Box$

\bigskip

We proceed to demonstrating the solvability in the sense of
sequences for our equation on the interval $I$ with periodic boundary
conditions.

\bigskip

\noindent
{\it Proof of Theorem 1.4.} According to our assumptions, we derive that
$$
|f(0)-f(2\pi)|\leq |f(0)-f_{m}(0)|+|f_{m}(2\pi)-f(2\pi)|\leq
2\|f_{m}-f\|_{C(I)}\to 0
$$
as $m\to \infty$. Thus, $f(0)=f(2\pi)$. By virtue of (\ref{f12}) for
$f_{m}(x), \ f(x)$ continuous on the interval $I$, we have
$f_{m}(x), \ f(x)\in L^{1}(I)\cap L^{2}(I), \ m\in {\mathbb N}$. First bound in
(\ref{f12}) also implies that
\begin{equation}
\label{fmfl1}
\|f_{m}(x)-f(x)\|_{L^{1}(I)}\leq 2\pi \|f_{m}(x)-f(x)\|_{C(I)}\to 0, \quad
m\to \infty.
\end{equation}
Therefore, $f_{m}(x)\to f(x)$ in $L^{1}(I)$ as $m\to \infty$.
Similarly, (\ref{f12}) gives us that
\begin{equation}
\label{fmfl2}
\|f_{m}(x)-f(x)\|_{L^{2}(I)}\leq \sqrt{2\pi} \|f_{m}(x)-f(x)\|_{C(I)}\to 0,
\quad m\to \infty.
\end{equation}
Hence, $f_{m}(x)\to f(x)$ in $L^{2}(I)$ as $m\to \infty$ as well.
Let us apply
the Fourier transform (\ref{fti}) to both sides of (\ref{eq5}). We
arrive at
\begin{equation}
\label{ufmn}
u_{m, n}=\frac{f_{m, n}}{n^{2}+n^{4}-a-ibn}, \quad m\in {\mathbb N}, \quad n\in
{\mathbb Z}.
\end{equation}
Let us consider the case a) of the theorem first. According to the part a) of
Theorem 1.3, equations (\ref{eq4}) and (\ref{eq5}) admit unique solutions
$u(x)\in H^{4}(I)$ and $u_{m}(x)\in H^{4}(I), \ m\in {\mathbb N}$ respectively.
Using (\ref{unfn}), (\ref{fmfl2}) and (\ref{ufmn}), we derive
$$
\|u_{m}-u\|_{L^{2}(I)}^{2}=\sum_{n=-\infty}^{\infty}
\frac{|f_{m, n}-f_{n}|^{2}}{(n^{2}+n^{4}-a)^{2}+b^{2}n^{2}}\leq \frac{1}{C}
\|f_{m}-f\|_{L^{2}(I)}^{2}\to 0
$$
as $m\to \infty$.
Thus, $u_{m}(x)\to u(x)$ in $L^{2}(I)$ as $m\to \infty$. We demonstrate
that $u_{m}(x)$ tends to $u(x)$ in $H^{4}(I)$ as
$m\to \infty$. It follows from (\ref{eq4}) and (\ref{eq5}) that
$$
\Bigg\|\Bigg(-\frac{d^{2}}{dx^{2}}\Bigg)(u_{m}-u)+\frac{d^{4}}{dx^{4}}(u_{m}-u)-b\frac{d(u_{m}-u)}{dx}\Bigg\|_
{L^{2}(I)}\leq 
$$
$$
\|f_{m}-f\|_{L^{2}(I)}+a\|u_{m}-u\|_{L^{2}(I)}.
$$
Let us recall (\ref{fmfl2}). Hence,
the right side of this inequality converges to zero as $m\to \infty$.
By means of the Fourier transform (\ref{fti}), we derive
$$
\sum_{n=-\infty}^{\infty}n^{8}|u_{m, n}-u_{n}|^{2}\to 0, \quad m\to \infty.
$$
Therefore,
$\displaystyle{\frac{d^{4}u_{m}}{dx^{4}}\to\frac{d^{4}u}{dx^{4}}}$ in
$L^{2}(I)$ as $m\to \infty$. This means that
$u_{m}(x)\to u(x)$ in $H^{4}(I)$ as $m\to \infty$ as well.

We conclude the work by discussing the situation when
the parameter $a$ is trivial. By virtue of (\ref{or4}) along with (\ref{fmfl1}),
we have
$$
|(f(x), 1)_{L^{2}(I)}|=|(f(x)-f_{m}(x), 1)_{L^{2}(I)}|\leq \|f_{m}-f\|_{L^{1}(I)}\to 0,
\quad m\to \infty,
$$
such that the limiting orthogonality relation
\begin{equation}
\label{oril}
(f(x), 1)_{L^{2}(I)}=0
\end{equation}
holds. Let us recall the part b) of Theorem 1.3 above. Thus, problems
(\ref{eq4}) and (\ref{eq5}) possess unique solutions $u(x)\in H_{0}^{4}(I)$
and $u_{m}(x)\in H_{0}^{4}(I), \ m\in {\mathbb N}$ respectively when $a$ vanishes.
Formulas (\ref{unfn}) and (\ref{ufmn}) give us that
\begin{equation}
\label{umnfmn}
u_{m, n}-u_{n}=\frac{f_{m, n}-f_{n}}{n^{2}+n^{4}-ibn}, \quad m\in {\mathbb N}, \quad n\in
{\mathbb Z}.
\end{equation}
Orthogonality conditions (\ref{oril}) and (\ref{or4}) imply
$$
f_{0}=0, \quad f_{m, 0}=0, \quad m\in {\mathbb N}.
$$
We obtain the estimate from above on the norm as
$$
\|u_{m}-u\|_{L^{2}(I)}=
\sqrt{\sum_{n=-\infty, \ n\neq 0}^{\infty}\frac{|f_{m, n}-f_{n}|^{2}}{(n^{2}+n^{4})^{2}+b^{2}n^{2}}}\leq
\frac{\|f_{m}-f\|_{L^{2}(I)}}{|b|}\to 0
$$
as $m\to \infty$ due to (\ref{fmfl2}). Therefore, $u_{m}(x)\to u(x)$ in $L^{2}(I)$ as $m\to \infty$.
Then $u_{m}(x)\to u(x)$ in $H_{0}^{4}(I)$ as $m\to \infty$ by
virtue of the reasoning analogous to the one above in the proof of the
case a) of our theorem.  

\hspace{12cm} $\Box$

\bigskip


\section{Acknowledgement}

\noindent
V.V. is grateful to Israel Michael Sigal for the
partial support by the NSERC grant.

\bigskip


\end{document}